\documentclass[a4paper,12pt]{amsart}
\usepackage[utf8x]{inputenc}
\usepackage{amsmath} 
\usepackage{amsthm}
\usepackage{amsfonts}
\usepackage{amssymb}
\usepackage{amscd}
\usepackage{mathrsfs}
\usepackage{graphicx}
\usepackage{longtable}
\usepackage{verbatim}
\usepackage[all]{xy}
\usepackage{hyperref}

\theoremstyle{plain}
\newtheorem{theorem}{Theorem}

\newtheorem{example}{Example} 
\newtheorem{question}{Question} 
\newtheorem{lemma}{Lemma}[section]

\newtheorem{proposition}[lemma]{Proposition}

\def\GL{\operatorname{GL}}

\def\SL{\operatorname{SL}}

\def\U{\operatorname{U}}

\def\Ad{\operatorname{Ad}}

\def\opd{\operatorname{d}}

\def\der{\operatorname{der}}
\def\det{\operatorname{det}}

\def\diag{\operatorname{diag}}

\def\Hom{\operatorname{Hom}}

\def\Ind{\operatorname{Ind}}

\def\nd{\operatorname{nd}}

\def\rank{\operatorname{rank}}

\def\span{\operatorname{span}}

\newcommand{\frh}{\mathfrak{h}}

\begin{document}

\title{A strong multiplicity one theorem for dimension data}
\date{February 24, 2022}
\thanks{}

\author{Jun Yu}
\address{Jun Yu, Beijing International Center for Mathematical Research and School of Mathematical Sciences, Peking
University, No. 5 Yiheyuan Road, Beijing 100871, China.}
\email{junyu@bicmr.pku.edu.cn}

\begin{abstract}
We call the dimension data $\mathscr{D}_{H_{1}}$ and $\mathscr{D}_{H_{2}}$ of two closed subgroups $H_{1}$ and $H_{2}$
of a given compact Lie group $G$ {\it almost equal} if $\mathscr{D}_{H_{1}}(\rho)=\mathscr{D}_{H_{2}}(\rho)$ for all
but finitely many irreducible complex linear representations $\rho$ of $G$ up to equivalence. When $G$ is connected,
we show that: if $\mathscr{D}_{H_{1}}$ and $\mathscr{D}_{H_{2}}$ are almost equal, then they are equal. When $G$ is
non-connected, $G^{0}\subset H_{1}\cap H_{2}$ is a trivial sufficient condition for $\mathscr{D}_{H_{1}}$ and
$\mathscr{D}_{H_{2}}$ to be almost equal. We show strong relations between $H_{1}$ and $H_{2}$ if $\mathscr{D}_{H_{1}}$
and $\mathscr{D}_{H_{2}}$ are almost equal but non-equal. We also construct an example to indicate that
$G^{0}\subset H_{1}\cap H_{2}$ is not a necessary condition.
\end{abstract}

\maketitle

\setcounter{tocdepth}{1}

\tableofcontents

\section*{Introduction}

Let $G$ be a given compact Lie group. Write $\hat{G}$ for the set of isomorphism classes of irreducible complex
linear representations of $G$, which is a countable set. The {\it dimension datum} of a closed subgroup $H$ is
defined by \[\mathscr{D}_{H}: \hat{G}\rightarrow\mathbb{Z},\quad\rho\mapsto\dim\rho^{H}.\] The dimension datum
was first studied by Larsen and Pink in their pioneering work \cite{Larsen-Pink}, with the motivation of helping
determine monodromy groups of $\ell$-adic Galois representations. In the beginning of the 21st century, Langlands
launched a program of ``beyond endoscopy", where he used dimension datum as a key ingredient in his stable trace
formula approach to showing general functoriality. Since then dimension datum catches more attention in the
mathematical community. Besides number theory and automorphic form theory, dimension datum also has applications
in differential geometry, e.g, it is used to construct non-isometric isospectral closed Riemannian manifolds
(\cite{An-Yu-Yu},\cite{Sutton}). In \cite{Yu-dimension}, we classified connected closed subgroups of a given
compact Lie group with the same dimension datum and we found generators of linear relations among dimension data
of connected closed subgroups.

We call the dimension data $\mathscr{D}_{H_{1}}$ and $\mathscr{D}_{H_{2}}$ of two closed subgroups $H_{1}$ and
$H_{2}$ of $G$ {\it almost equal} if $\mathscr{D}_{H_{1}}(\rho)=\mathscr{D}_{H_{2}}(\rho)$ for all but finitely
many $\rho\in\hat{G}$. In this paper we investigate the following problem: if $\mathscr{D}_{H_{1}}$ and
$\mathscr{D}_{H_{2}}$ are almost equal, are they actually equal? This might be regarded an analogue of the
strong multiplicity one theorem in automorphic form theory, which says that if two cuspidal automorphic
representations of $\GL(\mathbb{A}_{F})$ (for $F$ a number field) have the same Satake parameters at all but
finitely many places, then they are equal. When $G$ is connected, we show the following Theorem \ref{T1}.

\begin{theorem}\label{T1}
Assume that $G$ is connected. If the dimension data $\mathscr{D}_{H_{1}}$ and $\mathscr{D}_{H_{2}}$ of two closed
subgroups $H_{1}$ and $H_{2}$ of $G$ are almost equal, then they equal.
\end{theorem}

When $G$ is non-connected, $G^{0}\subset H_{1}\cap H_{2}$ is a trivial sufficient condition for $\mathscr{D}_{H_{1}}$
and $\mathscr{D}_{H_{2}}$ to be almost equal. In this case we show the following Theorem \ref{T2}. We also construct
an example which indicates that $G^{0}\subset H_{1}\cap H_{2}$ is not a necessary condition for $\mathscr{D}_{H_{1}}$
and $\mathscr{D}_{H_{2}}$ to be almost equal. 

\begin{theorem}\label{T2}
Assume that $Z(G^{0})^{0}\subset Z(G)$ or $G^0$ is semisimple. If $\mathscr{D}_{H_{1}}$ and $\mathscr{D}_{H_2}$ are
almost equal but non-equal, then $G^{0}\subset H_{1}\cap H_{2}$.
\end{theorem}

Similar questions are studied in the literature. For example, in \cite{Bhagwat-Rajan} it is shown a ``spectral analog
of strong multiplicity one theorem" for $L^{2}(\Gamma_{j}\backslash G)$ ($j=1,2$) with $G$ a real linear reductive
group and $\Gamma_{j}$ uniform lattices; in \cite{Lauret-Miatello} and \cite{Lauret-Miatello2}, Theorem \ref{T1} is
shown in the case that $H_{1},H_{2}$ are finite subgroups. The question studied in this paper might also be compared
to the following problem in isospectral geometry (\cite{Gordon}): if the Laplacian spectra of two closed Riemannian
manifolds have the same multiplies for all but finitely many eigenvalues, are they actually equal? To the
author's knowledge there is neither positive answer nor counter-example to this problem.

\smallskip

\noindent{\it Notation and conventions.} For a root system $\Phi$, let $\Phi^{\nd}$ denote the set of nondivisible
roots in $\Phi$, i.e., those roots $\alpha\in\Phi$ such that $\frac{1}{2}\alpha\not\in\Phi$.

Let $H$ be a closed abelian group. Write $X^{\ast}(H)=\Hom(H,\U(1))$ for the group of linear characters of $H$.
For a linear character $\lambda\in X^{\ast}(H)$, write $[\lambda]: H\rightarrow\U(1)\subset\mathbb{C}^{\times}$
for the corresponding evaluation function. Write $\mathbf{1}_{H}$ for the trivial character defined by 
\[[\mathbf{1}_{H}](x)=1,\ \forall x\in H.\] For two linear characters $\lambda,\mu\in X^{\ast}(H)$, define a 
linear character $\lambda+\mu$ by \[[\lambda+\mu](x)=[\lambda](x)[\mu](x),\ \forall x\in H.\] Write 
\[p: X^{\ast}(H)\rightarrow X^{\ast}(H^{0})\] be the projection defined by \[p(\lambda)=\lambda|_{H^{0}},
\ \forall\lambda\in X^{\ast}(H).\] 

Let $G$ be a compact Lie group. Let $\opd_{G}x$ be the unique normalized Haar measure on $G$. We call a conjugation
invariant function $\chi$ on $G$ a {\it virtual character} if $\span_{\mathbb{C}}\{g\cdot\chi\}$ is finite-dimensional,
where \[(g\cdot\chi)(g')=\chi(g^{-1}g'),\ \forall g'\in G.\] By representation
theory, any non-zero virtual character $\chi$ on $G$ is a finite sum \[\chi=\sum_{1\leq i\leq m}a_{i}\chi_{\rho_{i}}\]
where $\{\rho_{1},\dots,\rho_{m}\}$ is a finite set of pairwise non-isomorphic irreducible complex linear
representations of $G$ and $a_{1},\dots,a_{m}\in\mathbb{C}$ are non-zero constants.

\smallskip

\noindent{\it Acknowledgement.} This paper is an outcome of discussions between Yoshiki Oshima and the author.
In particular, Lemma \ref{L:data-measure} is due to him. I would like to thank Yoshiki Oshima for inspirational
discussions. I want to thank Emilio Lauret for helpful communications. This research is partially supported by the
NSFC Grant 11971036.

\section{Weyl integration formula and twining character formula}\label{S:Weyl}

In this section we recall a few results concerning non-connected compact Lie groups: parametrizing conjugacy classes,
Weyl density function and twining character formula. All results are from \cite{BtD}, \cite{Wendt}, 
\cite[\S 2,\S3]{Yu-dimension2}.

Let $H$ be a compact Lie group. We call a closed abelian subgroup $\tilde{S}$ of $H$ a {\it generalized Cartan subgroup}
if $\tilde{S}$ contains a dense cyclic subgroup and $\tilde{S}^{0}=Z_{H}(\tilde{S})^{0}$. We call a closed commutative
connected subset $S$ of $H$ a {\it maximal commutative connected subset} if $s^{-1}S=Z_{H}(S)^{0}$ for an (and hence
any) element $x\in S$. These objects are closely related: if $S$ is a maximal commutative connected subset, then
$\tilde{S}:=\langle S\rangle$ is a generalized Cartan subgroup; if $\tilde{S}$ is a generalized Cartan subgroup and
$x\in\tilde{S}$ generates $\tilde{S}/\tilde{S}^{0}$, then $S=x\tilde{S}^{0}$ is a maximal commutative connected subset.
Let $S$ be a maximal commutative connected subset of $H$. Set \[W(H,S)=N_{H^{0}}(S)/x^{-1}S\ (\forall x\in S),\] and
call it the {\it Weyl group} of $H$ with respect to $S$. Choose a maximal commutative connected subset in each connected
component of $H$. We obtain a set $\{S_{1},\dots,S_{k}\}$ ($k=|H/H^0|$). For each $i$ ($1\leq i\leq k$), write
$\tilde{S}_{i}=\langle S_{i}\rangle$. The following is a combination of some results in \cite[Chapter IV]{BtD}, which
parametrizes $H^{0}$ conjugacy classes in $H$.

\begin{proposition}[\cite{BtD}]\label{P:conjugacy}
Every element of $H$ is $H^{0}$ conjugate to an element in some $S_{i}$ ($1\leq i\leq k$). For each $i$, two elements
$x,x'\in S_{i}$ are $H^{0}$ conjugate if and only if they are in the same $W(H,S_{i})$ orbit.
\end{proposition} 

For each $i$ ($1\leq i\leq k$), put \begin{equation}\label{Eq:density}D_{H}(x)=\frac{1}{|W(H,S_{i})|}
\det(1-\Ad(x))_{\frh/\tilde{\mathfrak{s}}_{i}}|\ (\forall x\in S_{i}).\end{equation} We call $D_{H}(x)$ the
{\it Weyl density function} of $H$. Write $\opd x$ for the unique right $\tilde{S}_{i}^{0}$ invariant measure on
$S_{i}$ of volume $1$.

\begin{lemma}[\cite{Wendt}]\label{L:Weyl2}
Let $\chi$ be a conjugation invariant function on $H$. For any $H^{0}$ conjugation invariant continuous function $f$
on $H$, we have \[\int_{H}f(x)\chi(x)\opd_{H}x=\frac{1}{|H/H^{0}|}
\sum_{1\leq i\leq m}\int_{S_{i}}f(y)\chi(y)D_{H}(y)\opd y.\]
\end{lemma} 

Write $\tilde{R}(H,\tilde{S}_{i})$ for the set of non-trivial linear characters $\alpha\in X^{\ast}(\tilde{S}_{i})$ 
such that \[\mathfrak{h}_{\alpha}:=\{Y\in\mathfrak{h}:\Ad(x)Y=\alpha(x)Y,\ \forall x\in\tilde{S}_{i}\}\neq 0.\]
We call $\alpha\in\tilde{R}(H,\tilde{S}_{i})$ an infinite root if $\alpha|_{\tilde{S}_{i}^{0}}\neq\mathbf{1}$.
Let $R(H,\tilde{S}_{i})$ be the set of infinite roots in $\tilde{R}(H,\tilde{S}_{i})$ and we call it the
{\it quasi root system} of $H$ with respect to $\tilde{S}_{i}$. Put \[\bar{R}(H,\tilde{S}_{i})=
\{\alpha|_{\tilde{S}_{i}^{0}}:\alpha\in R(H,\tilde{S}_{i})\}\] and we call it the {\it restricted root system} of 
$H$ with respect to $\tilde{S}_{i}^{0}$, which is a root system in the lattice $X^{\ast}(\tilde{S}^{0})$ in the 
sense of \cite[Definition 2.2]{Yu-dimension}. 

Let $H_{i}$ be a closed subgroup of $H$ generated by $H^0$ and $S_{i}$. Put $m_{i}=|H_{i}/H_{i}^{0}|$. Define a
one-dimensional representation $\chi_{0,i}$ of $H_{i}$ by \[\chi_{0,i}(xy^{p})=e^{\frac{2p\pi\mathbf{i}}{m_{i}}},
\ \forall (x,y,p)\in H^{0}\times S\times\mathbb{Z}.\] Put $T_{i}=Z_{H^{0}}(\tilde{S}_{i}^{0})$ ($1\leq i\leq k$),
which is a maximal torus of $H^{0}$ (\cite{Yu-dimension2}). Write $R(H,T_{i})$ for the root system of $H$ with 
respect to the maximal torus $T_{i}$. Choose an element $s_{0,i}\in S_{i}$ such that the map \[p: \{\alpha\in R(H,\tilde{S}_{i}):\alpha(s_{0,i})=1\}\rightarrow\bar{R}(H,\tilde{S}_{i}),\quad\alpha\mapsto p(\alpha)=
\alpha|_{\tilde{S}_{i}^{0}}\] is a bijection, which is called a pinned element of $R(H,\tilde{S}_{i})$ in
\cite{Yu-dimension2}. Then, the adjoint action of $s_{0,i}$ on $\mathfrak{h}$ is a pinned automorphism
(\cite{Yu-dimension2}). Thus, the induced action of $s_{0,i}$ on $R(H,T_{i})$ stabilizes a positive system
$R(H,T_{i})^{+}$ of it and induces a permutation on the corresponding simple system. Let $\lambda$ be an
$s_{0,i}$ invariant integer dominant weight of $T_{i}$ with respect to $R(H,T_{i})^{+}$. Define a linear
character $\mu\in X^{\ast}(\tilde{S}_{i})$ by \[\mu(ts_{0,i}^{j})=\lambda(t),
\ (\forall (t,j)\in\tilde{S}_{i}^{0}\times\mathbb{Z}).\] Put \[R_{s_{0,i}}=\{\alpha\in R(H,\tilde{S}_{i}):
\alpha(s_{0,i})=1\}.\] Then, it is a root system and $p: R_{s_{0,i}}\rightarrow\bar{R}(H,\tilde{S}_{i})$ is 
an isomorphism. Choose a positive system $R_{s_{0,i}}^{+}$ of $R_{s_{0,i}}$. Put \[\delta_{i}=\frac{1}{2}
\sum_{\alpha\in R_{s_{0,i}}^{+}}m_{\alpha}\alpha\] and \begin{equation}\label{Eq:Ai-mu}A_{i}(\mu)=
\sum_{w\in W_{R_{s_{0,i}}}}\epsilon(w)[\mu+\delta_{i}-w\delta_{i}].\end{equation} There is a unique smooth 
function $\chi_{\mu}$ on $S_{i}$ determined by \begin{equation}\label{Eq:chi-mu}\chi_{\mu}
\prod_{\alpha\in {R_{s_{0,i}}}}(1-[-m_{\alpha}\alpha])=\sum_{w\in W_{R_{s_{0,i}}}}\epsilon(w)
[w\mu+w\delta_{i}-\delta_{i}].\end{equation} 

\begin{lemma}[\cite{Yu-dimension2}]\label{L:Weyl4}
Let $\lambda$ be an $s_{0,i}$ invariant integer dominant weight of $T_{i}$ with respect to $R(H,T_{i})^{+}$.
Then \[\chi_{\mu}(x)D_{H}(x)=\frac{1}{|W_{R_{s_{0,i}}}|}\sum_{w\in W_{R_{s_{0,i}}}}w\cdot A_{i}(\mu).\]
\end{lemma}

The following lemma follows from the so-called twining
character formula. The twining character formula is first shown by Jantzen (\cite{Jantzen}) and is proved by many
different methods in the literature. Interested readers could consult \cite{Kumar-Lusztig-Prasad}, \cite{Wendt}
or \cite{Yu-dimension2} for the statement of the twining character formula and the proof.

\begin{lemma}[\cite{Yu-dimension2}]\label{L:twining}
Let $\lambda$ be an $s_{0,i}$ invariant integer dominant weight of $T_{i}$ with respect to $R(H,T_{i})^{+}$. Then
there exists an irreducible complex linear representation $\sigma_{\lambda}$ of $H_{i}$ such that
\[\chi_{\sigma_{\lambda}}|_{S_{i}}=\eta\chi_{\mu}\] for some constant $\eta\in\mathbb{C}$ with $|\eta|=1$.

Let $\sigma$ be an irreducible complex linear representation of $H_{i}$. If $\sigma|_{H^{0}}$ has a dominant weight
which is not $s_{0,i}$ invariant, then $\chi_{\sigma}|_{S_{i}}=0$; if $\lambda$ is a dominant weight of
$\sigma|_{H^{0}}$, then $\sigma\cong\sigma_{\lambda}\otimes\chi_{0,i}^{j}$ for some $j\in\mathbb{Z}$.
\end{lemma}

\section{Sato-Tate measure associated to a closed subgroup}\label{S:STmeasure}

Let $G$ be a given compact Lie group. Endow $G$ with a biinvariant Riemannian metric induced by a conjugate invariant
inner product $(\cdot,\cdot)$ on $\mathfrak{g}_0$. Write $G^{\sharp}$ for the set of conjugacy classes in $G$, which
is a disjoint union of finitely many sets $S'_{i}/W'_{i}$ ($0\leq i\leq l$) with each $S'_{i}$ a maximal commutative
connected subset of $G$ and \[W'_{i}=N_{G}(S'_{i})/Z_{G}(S'_{i}).\] Here we choose $S'_{0}$ to be a maximal torus of
$G^{0}$. Let $2\delta'_{j}\in X^{\ast}(\tilde{S}'_{j})$ be the linear character associated to $G$ and
$\tilde{S}'_{j}$ analogous to that of $2\delta_{i}\in X^{\ast}(\tilde{S}_{i})$ associated to $H$ and
$\tilde{S}_{i}$ as in Section \ref{S:Weyl}. 

For a closed subgroup $H$ of $G$, define a measure $\mu_{H}$ on $G^{\sharp}$ by \[\mu_{H}(f)=
\int_{H}f(x)\opd_{H} x\ (\forall f\in C(G)^{G}),\] which is called the {\it Sato-Tate measure of $H$}. As in
Section \ref{S:Weyl}, we choose a maximal commutative connected subset in each connected component of $H$,
and obtain a set $\{S_{1},\dots,S_{k}\}$ ($k=|H/H^0|$). For each $i$ ($1\leq i\leq k$), let $\tilde{S}_{i}$
be the closed subgroup of $H$ generated by $S_{i}$, which is a generalized Cartan subgroup of $H$. For each
$i$, there exists a unique $j$ ($0\leq j\leq l$) such that $S_{i}$ is $G$-conjugate to a subset of
$S'_{j}$. By Lemma \ref{L:Weyl2}, the restriction of the measure $\mu_{H}$ to $S_{i}$ is a multiple of
$D_{H}(x))|_{S_{i}}$.

\begin{lemma}\label{L:mu-H}
We have the following assertions. \begin{enumerate}
\item[(i)]We have $\dim S_{i}\leq\dim S'_{j}$.
\item[(ii)]If $\dim S_{i}=\dim S'_{j}$, then the leading term appearing in $D_{H}(x)|_{S_{i}}$ has length
$\leq |2\delta'_{j}|$ and the equality holds only when $(Z(G^{0})^{x})^{0}(G^{0})_{\der}$ ($\forall x\in S'_{j})$.
\end{enumerate} 
\end{lemma}

\begin{proof}
(i)Since $S_{i}$ is $G$-conjugate to a subset of $S'_{j}$, we get $\dim S_{i}\leq\dim S'_{j}$.

(ii)If $\dim S_{i}=\dim S'_{j}$, then $S_{i}$ is equal to a $G$-conjugate of $ S'_{j}$. Without loss of generality
we assume that $S_{i}= S'_{j}$. Then, $\tilde{S}_{i}^{0}=\tilde{S}_{j}^{'0}$. Replacing $H$ by a conjugate one if 
necessary, we may assume that $(2\delta_{i})|_{\tilde{S}_{i}^{0}}$ is dominant with respect to a positive system 
$\bar{R}(G,\tilde{S}'_{j})^{+}$ of $\bar{R}(G,\tilde{S}'_{j})$. Since $((2\delta_{i})|_{\tilde{S}_{i}^{0}},
\bar{\alpha})>0$ for any root $\bar{\alpha}\in\bar{R}(H,\tilde{S}_{i})^{+}$, we have 
$\bar{R}(H,\tilde{S}_{i})^{+}\subset\bar{R}(G,\tilde{S}'_{j})^{+}$. Then, \[|2\delta'_{j}|^{2}-|2\delta_{i}|^{2}
=(2\delta'_{j}+2\delta_{i},2\delta'_{j}-2\delta_{i})\geq 0.\] Thus, $|2\delta_{i}|\leq|2\delta'_{j}|$. If the 
equality holds, then $2\delta'_{j}-2\delta_{i}=0$. Thus, \[\bar{R}(H,\tilde{S}_{i})^{+}=
\bar{R}(G,\tilde{S}'_{j})^{+},\] which indicates that $(Z(G^{0})^{x})^{0}(G^{0})_{\der}\subset H$ 
($\forall x\in S'_{j})$. 
\end{proof}

When $G$ is connected, we have $l=0$ and $ S'_0$ is a maximal torus of $G$. The following Lemma
\ref{L:connected-H} strengthens Lemma \ref{L:mu-H} when $G$ is connected.

\begin{lemma}\label{L:connected-H}
Let $G$ be a connected compact Lie group and $H$ be a closed subgroup. Then we have the following assertions.
\begin{enumerate}
\item[(i)]Except for probably maximal tori of $H^0$, all other maximal commutative connected subsets of $H$
have dimension strictly lower than $\dim S'_{0}$.
\item[(ii)]If a maximal torus $S_{i}$ of $H^0$ has dimension equal to $\dim S'_{0}$, then the dominant term of
the restriction of the measure $\mu_{H}$ to $S_{i}$ has length $\leq|2\delta_{0}|$ and the equality holds only
when $H=G$.
\end{enumerate}
\end{lemma}

\begin{proof}
(i)Let $S$ be a maximal commutative connected subset of $H$, not contained in $H^{0}$. Since $G$ is connected,
there must exist a maximal torus $T$ of $G$ containing $S$. Then, $S=xT'$ ($\forall x\in S$) with $T'$ a
sub-torus of $T$. Since $1\not\in S$, we have $T'\neq T$. Thus, \[\dim S=\dim T'<\dim T=\dim S'_{0}.\]

(ii)This follows from Lemma \ref{L:mu-H} (ii) directly.
\end{proof}

\section{Measure associated to a virtual character}\label{S:measure}

For a virtual character $\chi\in\mathbb{C}[\hat{G}]$, define a measure $\mu(\chi)$ on $G^{\sharp}$ by \[\mu(\chi)(f)
=\int_{G}f(x)\chi(x)\opd_{G} x\ (\forall f\in C(G)^{G}),\] where $C(G)^{G}$ means the space of conjugate invariant
continuous functions on $G$. Let $\pi_{0}(G)=G/G^{0}$ be the component group of $G$. By Lemma \ref{L:Weyl2}, the
restriction of the measure $\mu(\chi)$ to $S'_{j}$ is given by a multiple of $\chi(x)D_{G}(x)|_{S'_{j}}$. In the
following Lemma \ref{L:mu-chi} we show a lower bound for the length of the leading term appearing in
$(\chi(x)D_{G}(x))|_{S'_{j}}$ when $\chi|_{S'_{j}}\neq 0$. 

\begin{lemma}\label{L:mu-chi}
We have the following assertions. \begin{enumerate}
\item[(i)]For any virtual character $\chi$ of $G$, $\chi|_{S'_{j}}$ is a finite sum of terms $a_{\mu}\chi_{\mu}$,
where $a_{\mu}\in\mathbb{C}$ and $\chi_{\mu}$ is defined as in Section \ref{S:Weyl} (for the group
$\langle G^{0},S'_{j}\rangle$).
\item[(ii)]If $\chi|_{S'_{j}}\neq 0$, then the leading term appearing in $\chi(x)D_{G}(x)|_{S'_{j}}$ has
length $\geq|2\delta'_{j}|$. Moreover, if the equality holds, then $\chi|_{S'_{j}}$ is a nonzero constant.
\end{enumerate}
\end{lemma}

\begin{proof}
(i)By decomposing $\sigma|_{G_{j}}$ for each irreducible representation $\sigma$ of $G$ appearing in the virtual
character $\chi$, we write $\chi|_{S'_{j}}$ as a finite sum of terms $a_{\mu}\chi_{\mu}$ by Lemma \ref{L:twining}.

(ii)By Lemma \ref{L:Weyl4}, the leading term appearing in $\chi(x)D_{G}(x)|_{S'_{j}}$ is the longest $\mu+
2\delta'_{j}$ among those dominant integral weights $\mu$ such that $\chi_{\mu}$ appears in $\chi|_{S'_{j}}$.
Hence, the leading term appearing in $\chi(x)D_{G}(x)|_{S'_{j}}$ has length $\geq|2\delta'_{j}|$. If the
equality holds, one must have all $\mu=0$. Then, $\chi|_{S'_{j}}$ is a nonzero constant.
\end{proof}

When $G$ is connected, Lemma \ref{L:mu-chi} specializes to the following Lemma \ref{L:connected-chi}.

\begin{lemma}\label{L:connected-chi}
Let $G$ be a connected compact Lie group. Then we have the following assertions. \begin{enumerate}
\item[(i)]For any nonzero virtual character $\chi$ of $G$, the measure $\mu(\chi)$ is given by integration against
the function \[\frac{1}{|W(G,S'_0)|}\chi(x)D_{G}(x)\opd x\] on $S'_0$.
\item[(ii)]The leading term of $\chi(x)D_{G}(x)$ has length $\geq|2\delta'_{0}|$ and the equality holds only when
$\chi$ is a nonzero constant function. 
\end{enumerate}
\end{lemma}

\section{Closed subgroups with almost equal dimension data}\label{S:strongMultiplicity}

Let $H_{1},H_{2}$ be two closed subgroups of $G$. If $\mathscr{D}_{H_1}$ and $\mathscr{D}_{H_2}$ are almost equal,
then \[X=\{\rho\in\hat{G}:\mathscr{D}_{H_1}(\rho)\neq\mathscr{D}_{H_2}(\rho)\}\] is a finite set. For each
$\rho\in X$, write \[a_{\rho}=\mathscr{D}_{H_1}(\rho)-\mathscr{D}_{H_2}(\rho).\] Put \[\chi:=\sum_{\rho\in X}
a_{\rho}\chi_{\rho^{\ast}}=\sum_{\rho\in X}a_{\rho}\overline{\chi_{\rho}}.\] We call $X$ the exceptional support,
and call $\chi$ the exceptional character. Apparently, $\mathscr{D}_{H_1}=\mathscr{D}_{H_2}$ if and only if
$X=\emptyset$, or equivalently $\chi=0$.

\begin{lemma}\label{L:data-measure}
For $\mathscr{D}_{H_1}$ and $\mathscr{D}_{H_2}$ to be almost equal, it is necessary and sufficient that there is a virtual
character $\chi$ such that $\mu_{H_1}-\mu_{H_2}=\mu(\chi).$
\end{lemma}

\begin{proof}
Let $X$ be a finite subset of $\hat{G}$. Write \[\chi=\sum_{\rho\in X}a_{\rho}\chi_{\rho^{\ast}}.\] By the orthogonality
theorem, for each representation $\rho\in\hat{G}$, one has \[\mu(\chi)(\chi_{\rho})=\left\{\begin{array}{cc}a_{\rho}\quad
if\ \rho\in X\\0\quad if\ \rho\not\in X.\\\end{array}\right.\] On the other hand, \[\mu_{H_1}(\chi_{\rho})-
\mu_{H_2}(\chi_{\rho})=\int_{H_{1}}\chi_{\rho}(x)\opd_{H_{1}}x-\int_{H_{2}}\chi_{\rho}(x)\opd_{H_{2}}x=
\mathscr{D}_{H_{1}}(\rho)-\mathscr{D}_{H_{2}}(\rho).\] Therefore, \[\mathscr{D}_{H_{1}}(\rho)-\mathscr{D}_{H_{2}}(\rho)
=\left\{\begin{array}{cc}0\quad if\ \rho\in\hat{G}-X\\a_{\rho}\quad if\ \rho\in X\qquad\\\end{array}\right.\] is equivalent
to \[\mu_{H_1}-\mu_{H_2}=\mu(\chi).\]
\end{proof}

With Lemma \ref{L:data-measure}, we transfer the question of whether $\mathscr{D}_{H_1}$ and $\mathscr{D}_{H_2}$ are
almost equal to the question of whether $\mu_{H_1}-\mu_{H_2}=\mu(\chi)$ for a virtual character $\chi$. When $G$ is
connected, an affirmative answer is given in Theorem \ref{T1}.

\begin{proof}[Proof of Theorem \ref{T1}]
Suppose that $\mathscr{D}_{H_1}$ and $\mathscr{D}_{H_2}$ are almost equal, but non-equal. By Lemma \ref{L:data-measure},
there is a nonzero virtual character $\chi$ such that $\mu_{H_1}-\mu_{H_2}=\mu(\chi)$. We compare the leading term in
the $\rank G$-dimensional component of the supports of $\mu_{H_1}-\mu_{H_2}$ and $\mu(\chi)$. By Lemma
\ref{L:connected-chi} and Lemma \ref{L:connected-H} we must have $H_1=G$ or $H_{2}=G$. Without loss of generality, we
assume that $H_2=G$. Then, \[\mu(\chi+1)=\mu_{H_2}+\mu(\chi)=\mu_{H_{1}}\neq 0.\] Comparing the leading term in the
$\rank G$-dimensional component of the supports of $\mu_{H_{1}}$ and $\mu(\chi+1)$ again, we get $H_1=G$. This is in
contradiction with $\mathscr{D}_{H_1}\neq\mathscr{D}_{H_2}$.
\end{proof}

\begin{proposition}\label{P:strong2}
Let $H_{1},H_{2}$ be closed subgroups of $G$. Suppose that $\mu_{H_{1}}-\mu_{H_2}=\mu(\chi)$ for a non-zero virtual
character $\chi$ of $G$. Then, we have the following assertions. \begin{enumerate}
\item[(i)]Every irreducible representation $\rho$ appearing in $\chi$ is trivial on $G^{0}$.
\item[(ii)]$Z(G)^{0}(G^{0})_{\der}\subset H_{1}\cap H_{2}$.
\item[(iii)]If $\chi|_{G^{0}}\neq 0$, then $G^{0}\subset H_{1}\cap H_{2}$.
\item[(iv)]If $\chi|_{G^{0}}=0$, then $\mu_{H_{1}\cap G^{0}}=\mu_{H_{2}\cap G^{0}}$ and $|H_{1}/H_{1}\cap G^0|=|H_{2}/H_{2}
\cap G^0|$.
\end{enumerate}
\end{proposition}

\begin{proof}
(i)Take a connected component $g_{j}G^{0}$ such that $\chi|_{g_{j}G^{0}}\neq 0$. Take a maximal commutative connected
subset $S'_{j}$ in $g_{j}G^{0}$. Consider the restriction of $\mu_{H_{1}}-\mu_{H_{2}}$ and $\mu(\chi)$ to $S'_{j}$ and
the leading terms of their expressions. By Lemma \ref{L:mu-chi} and Lemma \ref{L:mu-H}, one sees that $\chi|_{g_{j}G^{0}}$
is a constant function and at least one of $H_{1}$ and $H_{2}$ contains $(Z(G^{0})^{g_{j}})^{0}(G^{0})_{\der}$. Letting
$j$ vary, we see that the restriction of $\chi$ to each connected component of $G$ is a constant function. Hence, $\chi$
is the pull-back of a virtual character of $\pi_{0}(G)$. Therefore, every representation $\rho$ appearing in $\chi$ is
trivial on $G^{0}$.

(ii)Without loss of generality we suppose that $Z(G)^{0}(G^{0})_{\der}\not\subset H_{2}$. Continue the argument in (i),
but now we further consider the coefficients of the leading terms of expressions of the restrictions of $\mu_{H_{1}}-
\mu_{H_{2}}$ and $\mu(\chi)$ to $S'_{j}$, we get $\chi(g_{j})>0$. Letting $g_{j}$ vary, we see that $\chi$
is a positive character. Then, $a_{\mathbf{1}_{G}}=(\mathbf{1}_{\tilde{G}},\chi)>0$, which is in contradiction with \[a_{\mathbf{1}_{G}}=\mathscr{D}_{H_1}(\mathbf{1}_{G})-\mathscr{D}_{H_2}(\mathbf{1}_{G})=0.\]

(iii)Suppose that $\chi|_{G^{0}}\neq 0$. Then, one shows $G^{0}\subset H_{1}\cap H_{2}$ similarly as in the proof of
Theorem \ref{T1}.

(iv) Suppose that $\chi|_{G^{0}}=0$. By considering the restriction of $\mu_{H_1},\mu_{H_{2}}$ to a maximal torus of
$G^{0}$, we get \[\frac{1}{|H_{1}/H_{1}\cap G^0|}\mu_{H_{1}\cap G^{0}}=\frac{1}{|H_{2}/H_{2}\cap G^0|}
\mu_{H_{2}\cap G^{0}}.\] Evaluating at the function $\mathbf{1}_{G^{0}}$, we get \[|H_{1}/H_{1}\cap G^0|=
|H_{2}/H_{2}\cap G^0|.\] Thus, \[\mu_{H_{1}\cap G^{0}}=\mu_{H_{2}\cap G^{0}}.\]
\end{proof} 

\begin{proof}[Proof of Theorem \ref{T2}]
By Lemma \ref{L:data-measure} and Proposition \ref{P:strong2} (ii), we get \[Z(G)^{0}(G^{0})_{\der}
\subset H^{1}\cap H^{2}.\] Since $Z(G^{0})^{0}\subset Z(G)$ or $G^0$ is semisimple, we have
\[Z(G)^{0}(G^{0})_{\der}=G^0.\] Then, $G^{0}\subset H_{1}\cap H_{2}$.
\end{proof}

On the other hand, suppose that $G^{0}\subset H^{1}\cap H^{2}$. Then $\mathscr{D}_{H_{1}}(\rho)=\mathscr{D}_{H_2}(\rho)=0$
for an irreducible representation $\rho$ whenever $\rho|_{G^{0}}$ is non-trivial. Particularly, $\mathscr{D}_{H_{1}}(\rho)-
\mathscr{D}_{H_2}(\rho)\neq 0$ only for irreducible representations which factor through $G/G^0$. Hence,
$\mathscr{D}_{H_{1}}$ and $\mathscr{D}_{H_2}$ are almost equal.

\begin{lemma}\label{L:strong4}
Let $n>1$ and $\Gamma$ be a non-trivial finite subgroup of $\SL(n,\mathbb{Z})$. Suppose that any $1\neq\gamma\in\Gamma$ has
no eigenvalue $1$ while acting on $\mathbb{C}^{n}$. Form \[G=\U(1)^{n}\rtimes\Gamma.\] Take two subgroups $H_{1},H_{2}$ of
$\Gamma$ with the same order. Consider them as finite subgroups of $G$. Then, $\mathscr{D}_{H_{1}}(\rho)-
\mathscr{D}_{H_{2}}(\rho)=0$ for any irreducible representation $\rho$ of $G$ whenever $\rho|_{G^{0}}$ is nontrivial.
\end{lemma}

\begin{proof}
Since any $1\neq\gamma\in\Gamma$ has eigenvalue $0$, elements in $gG^{0}$ consists of a single $G^{0}$ orbit for any $g\in
G-G^{0}$. Hence, $G^{\sharp}$ is the union of $G^{0}/\sim_{G}$ and finitely many points. Since $H_1$ and $H_{2}$ have the
same order, the restriction of the measure $\mu_{H_{1}}-\mu_{H_{2}}$ to $G^{0}/\sim_{G}$ is $0$. Take a $G^{0}$ linear
character $\sigma\subset\rho|_{G^{0}}$, which is nontrivial as it is assume that $\rho|_{G^{0}}$ is nontrivial. Then, by the
Frobenius reciprocity one has $\rho\subset\Ind_{G^0}^{G}\sigma$. Since any $1\neq\gamma$($\in\Gamma$) has no eigenvalue $1$,
one has $(\gamma_{2}^{-1}\gamma_{1})\cdot\sigma\neq\sigma$ for any pair of distinct $\gamma_{1},\gamma_{2}\in\Gamma$. Hence, $\gamma_{1}\cdot\sigma\neq\gamma_{2}\cdot\sigma$. Then, $\Ind_{G^0}^{G}\sigma$ is irreducible. Therefore, $\rho=
\Ind_{G^0}^{G}\sigma$. By this, $\chi_{\rho}$ is supported in $G^{0}$. Therefore, \[\mathscr{D}_{H_{1}}(\rho)-
\mathscr{D}_{H_{2}}(\rho)=(\mu_{H_{1}}-\mu_{H_{2}},\chi_{\rho})=0.\]
\end{proof}

\begin{example}\label{E:strong5}
Write $J_{m}=\left(\begin{array}{cc}0_{m}&I_{m}\\-I_{m}&0_{m}\end{array}\right).$ Let $n=4$. Take $$\Gamma\!=\!\langle\diag
\{J_1,-J_{1}\},J_{2}\rangle\cong Q_{4},$$ $$H_{1}=\langle\diag\{J_1,-J_{1}\}\rangle\cong\mathbb{Z}/4\mathbb{Z},$$ $$H_2=\langle J_2\rangle\cong\mathbb{Z}/4\mathbb{Z}.$$ Then, the finite group $\Gamma$ satisfies the condition in Lemma \ref{L:strong4} and
$|H_1|=|H_{2}|$. Hence, $\mathscr{D}_{H_{1}}$ and $\mathscr{D}_{H_{2}}$ are almost equal. On the other hand, since the supports
of $\mu_{H_1}$ and $\mu_{H_2}$ are non-equal, one has $\mathscr{D}_{H_1}\neq\mathscr{D}_{H_2}$.
\end{example}

\begin{question}\label{Q:strong5}
Let $G$ be a compact Lie group. Can one give a group-theoretical criterion for pairs of closed subgroups $H_{1},H_{2}$
such that \[\mu_{H_{1}}-\mu_{H_{2}}=\mu(\chi)\] with $\chi$ a non-zero virtual character?
\end{question}

With Proposition \ref{P:strong2}, by modulo $Z(G)^{0}(G^{0})_{\der}$ one reduces Question \ref{Q:strong5} to the case that $G^{0}$
is a torus and $Z(G)$ is a finite group. By Example \ref{E:strong5}, it is not necessarily that $G^{0}\subset H_{1}\cap H_{2}$
in Question \ref{Q:strong5}.

\end{document}